\documentclass[11pt]{article}
\usepackage{booktabs}
\usepackage{makecell}
\usepackage{booktabs}
\usepackage{siunitx}
\usepackage{array}
\usepackage{amsmath}
\usepackage{amssymb}
\usepackage{algorithm}
\usepackage{algpseudocode}
\usepackage{bbm}
\usepackage{graphicx} % Required for inserting images
\usepackage{subcaption} % for subfigures
\usepackage[colorlinks=true, linkcolor=blue, citecolor=blue, urlcolor=blue]{hyperref}

\usepackage[table]{xcolor}
\usepackage{slashbox}
\usepackage{setspace}
\usepackage[a4paper, margin=1in]{geometry}
\usepackage{multirow}
\usepackage{geometry}
\usepackage{amsthm}

\usepackage{graphicx}

\usepackage{amsmath}
\usepackage{amsfonts}
\usepackage{amssymb}
\usepackage{amsthm}

\usepackage{url}
\usepackage{fancyhdr}
\usepackage{indentfirst}
\usepackage{enumerate}
\usepackage{booktabs}
\usepackage{color}
\usepackage{threeparttable}
\usepackage{float}
\usepackage{comment}

\usepackage{algorithm}
\usepackage{algpseudocode}

\usepackage{tikz}
\usetikzlibrary{positioning,arrows.meta,fit}

\newcommand{\be}{\begin{eqnarray}}
\newcommand{\ee}{\end{eqnarray}}
\newcommand{\by}{\begin{eqnarray*}}
\newcommand{\ey}{\end{eqnarray*}} % <--- CORRECTED LINE HERE
\newcommand{\bn}{\begin{enumerate}}
\newcommand{\en}{\end{enumerate}}
\newcommand{\bi}{\begin{itemize}}
\newcommand{\ei}{\end{itemize}}

\renewcommand{\ge}{\geqslant}
\renewcommand{\le}{\leqslant}

\renewcommand{\epsilon}{\varepsilon}
\theoremstyle{plain}
\newtheorem{theorem}{Theorem}

\theoremstyle{definition}

\numberwithin{equation}{section} \numberwithin{theorem}{section}
\renewcommand{\cite}{\citet} % This line makes \cite work like \citet from natbib.

\newcommand{\ignore}[1]{}

\usepackage{newtxtext}
\usepackage{newtxmath}
\usepackage{bm}

\usepackage[margin=1in]{geometry}

\usepackage[colorlinks=true,citecolor=blue,linkcolor=blue]{hyperref}

\usepackage{fancyhdr}
\title{Sampler-Robust Optimization under Generative Models
}
\author{
\large Ziwei Zhang \hspace{3em} Jonathan Yu-Meng Li\\[1em]
\normalsize Telfer School of Management\\
University of Ottawa\\
\normalsize\texttt{zzhan073@uottawa.ca, jonathan.li@telfer.uottawa.ca}
}

\usepackage{natbib}
\begin{document}
\maketitle
\begin{abstract}
Modern stochastic optimization pipelines increasingly rely on learned generative models to represent uncertainty, while downstream decisions are evaluated almost entirely through Monte Carlo scenarios. This shifts the operational object of uncertainty from an explicit probability law to the sampler induced by the learned generator. Reliability therefore depends on two errors: sampler misspecification and finite-simulation error. We propose \emph{Sampler-Robust Optimization} (SRO), which optimizes decisions against the worst-case sampler induced by perturbing the learned generator. This sampler-first formulation aligns with simulation-based decision pipelines and admits a sharpness-aware interpretation: it favors decisions whose performance is stable under generator perturbations, rather than merely under the nominal sampler. Under a coverage assumption, we show that the empirical worst-case objective provides a high-probability upper certificate for the true population objective, with finite-simulation error partially absorbed by the robustification used to guard against sampler misspecification. The framework accommodates generative models with or without explicit densities and admits efficient minimax procedures. Portfolio-optimization experiments show that SRO produces more stable decisions and improves out-of-sample performance under distribution shift.

\end{abstract}

\section{Introduction}

Generative models are increasingly used as uncertainty engines in modern decision pipelines. In applications such as portfolio allocation, supply-chain planning, energy systems, and reinforcement learning, conditional generative models are trained to simulate future outcomes from historical and contextual information \citep{ni2024sigwgan,li2025gmdm}. Once deployed, however, these models rarely enter downstream decision problems through direct manipulation of an explicit probability law. Instead, objectives and constraints are implemented almost entirely through Monte Carlo scenarios generated by the learned model.

This operational fact changes the object that truly enters the optimization problem. When a generative model is embedded in a simulation-based decision pipeline, the optimizer does not access an abstract distribution in closed form; it interacts with the \emph{sampler} induced by the learned generator. Candidate decisions are evaluated by drawing latent noise, passing it through the generator, and aggregating the resulting losses. Reliability therefore depends on two distinct layers. First, the learned sampler itself may be misspecified, in that the fitted generator may not coincide with the true data-generating mechanism. Second, even under a fixed sampler, decision quality is assessed from finitely many generated scenarios, creating simulation error.

Most existing robust formulations do not directly target this sampler-level uncertainty. Classical distributionally robust optimization constructs ambiguity sets around probability distributions \citep{ben2013robust,gao2016distributionally,esfahani2018data}. More recent approaches parameterize adversarial distributions through learned generators, including KL-constrained generative DRO formulations \citep{michel2021modeling} and generative ambiguity-set methods such as GAS-DRO \citep{wen2026gasdro}. Despite using generators, these formulations remain distribution-centered, relying on likelihood-based ingredients to define or compare adversarial laws. In many modern decision pipelines, however, the downstream optimizer never accesses such laws explicitly; it only queries a simulator. This distinction is especially pronounced in sampling-first settings such as simulation-based inference, where the generative system is primarily available as a black-box sampler \citep{cranmer2020frontier}. This motivates a different starting point: if the sampler is the operational object used for decision-making, then robustness should be imposed directly at the sampler level.

We pursue this idea through \emph{Sampler-Robust Optimization} (SRO), which optimizes decisions against the worst-case sampler induced by perturbing the learned generative model. SRO admits a sharpness-aware interpretation: rather than selecting a decision that performs well under a single learned sampler, it favors decisions whose performance remains stable under samplers induced by nearby generator perturbations. This connects SRO to perturbation-based methods such as Sharpness-Aware Minimization (SAM) \citep{foret2021sharpness}, but with a different purpose: decision-time robustness rather than training-time regularization.
We show that this robustification does more than hedge against sampler misspecification: under a coverage assumption, the empirical worst-case objective yields a high-probability upper certificate for the true population objective, with finite-simulation error partially absorbed by the robustification. The same mechanism that protects against misspecification therefore also improves reliability against Monte Carlo uncertainty. This reveals a perturbation-radius tradeoff: enlarging the neighborhood makes optimization more conservative, but can strengthen the resulting performance certificate. The framework applies broadly to generative models with or without explicit densities and admits efficient minimax procedures.

As a motivating application, we study portfolio optimization under learned multi-asset return generative models. Portfolio decisions are well known to be highly sensitive to misspecification of the return-generating mechanism, especially in high-dimensional and non-stationary environments \citep{tu2011markowitz,kan2007optimal}. Modern conditional generative models can represent rich cross-asset dependence, including through GAN-based architectures for financial time series \citep{mirza2014conditional,yoon2019time,vuletic2024fin}. Yet such models are often trained on limited historical trajectories and may produce unstable decisions under regime shift. More broadly, in complex non-stationary environments, even small perturbations to a learned generator can induce substantial changes in sampled outcomes, leading to fragile downstream decisions \citep{cont2001empirical}. SRO is designed precisely to address this vulnerability by robustifying decisions against perturbations of the learned sampler that generates future return scenarios.

\paragraph{Contributions.}
Our contributions are threefold.
\begin{itemize}
    \item We introduce a \emph{sampler-first} robust optimization framework for decision-making with learned generative models, in which robustness is imposed at the level of the induced sampler rather than through an explicit probability law.

    \item We show that sampler-level robustification naturally yields a finite-sample reliability guarantee that remains meaningful for large, highly parameterized generators. Under a coverage assumption, the empirical worst-case objective provides a high-probability upper certificate for the true population objective, with finite-simulation error partially absorbed by the robustification.

    \item We develop efficient minimax solution procedures and show, through portfolio-optimization experiments, that SRO improves decision stability and out-of-sample performance under distribution shift.
\end{itemize}

\section{Sampler-Robust Stochastic Optimization}
\label{sec:formulation}

We consider a context-dependent stochastic decision problem. Let $\mathcal X$ denote the context space, $\mathcal Y$ the outcome space, $\mathcal W \subseteq \mathbb R^d$ the feasible decision set, and
\[
f:\mathcal W \times \mathcal Y \to \mathbb R
\]
a measurable loss function that evaluates the performance of decision $\omega \in \mathcal W$ under outcome $y \in \mathcal Y$. Modern stochastic optimization pipelines are increasingly coupled with learned conditional generative models,
\[
P_\theta(\cdot \mid x),
\]
parameterized by $\theta \in \Theta \subseteq \mathbb R^p$, to represent the conditional distribution over outcomes in $\mathcal Y$ given contextual information $x \in \mathcal X$.

Under a nominal estimate $\hat\theta$, the population objective is
\begin{equation}
J(\omega;x)
:= \mathbb E_{y \sim P_{\hat\theta}(\cdot\mid x)}[f(\omega,y)],
\label{eq:nominal_population_objective}
\end{equation}
and the associated nominal decision rule is
\begin{equation}
\omega^*(x)\in \arg\min_{\omega\in\mathcal W} J(\omega;x).
\label{eq:nominal_decision_problem}
\end{equation}

Although \eqref{eq:nominal_population_objective} is written in terms of a probability law, it is typically evaluated in practice through Monte Carlo simulation. Thus, the learned generative model enters the decision problem only through the \emph{sampler} it induces. This motivates a formulation in which robustness is imposed directly at the sampler level, rather than through an abstract probability distribution.

\subsection{Sampler-first formulation}

Assume the conditional model admits a generative representation
\begin{equation}
y = G_\theta(z,x), \qquad z \sim \nu,
\label{eq:generative_representation}
\end{equation}
where $\nu$ is a fixed latent distribution on a latent space $\mathcal Z$, and
\[
G_\theta:\mathcal Z\times\mathcal X\to\mathcal Y
\]
is a measurable generator. Under \eqref{eq:generative_representation}, the population objective in \eqref{eq:nominal_population_objective} can be written equivalently as
\begin{equation}
J(\omega;x)
=
\mathbb E_{z\sim\nu}\!\left[f\bigl(\omega,G_{\hat\theta}(z,x)\bigr)\right].
\label{eq:nominal_population_sampler}
\end{equation}
Given i.i.d.\ latent draws $z_1,\ldots,z_N\sim \nu$, the corresponding empirical objective is
\begin{equation}
\hat J_N(\omega;x)
=
\frac{1}{N}\sum_{i=1}^N f\bigl(\omega,G_{\hat\theta}(z_i,x)\bigr).
\label{eq:nominal_empirical_sampler}
\end{equation}

Equations \eqref{eq:nominal_population_sampler}--\eqref{eq:nominal_empirical_sampler} make explicit that, once the model is queried only through simulation, the optimizer interacts not with a closed-form probability law but with the induced sampler $G_{\hat\theta}$. This immediately reveals two distinct sources of unreliability:
\begin{enumerate}
    \item \textbf{Sampler misspecification:} the learned parameter $\hat\theta$ may differ from the true parameter $\theta^\star$;
    \item \textbf{Simulation error:} the objective is evaluated using finitely many latent samples.
\end{enumerate}
Our aim is to address both within a single decision framework.

\subsection{Worst-case samplers}

To account for sampler misspecification, we begin with a perturbation neighborhood around the learned generator parameter,
\begin{equation}
\Theta_\rho(\hat\theta)
:=
\bigl\{\theta\in\Theta:\|\theta-\hat\theta\|\le \rho\bigr\},
\label{eq:theta_ball}
\end{equation}
where $\|\cdot\|$ is a chosen norm on $\mathbb R^p$ and $\rho>0$ controls the perturbation radius. For any decision $\omega$, we define the robust population objective
\begin{equation}
U(\omega; x)
:=
\sup_{\theta \in \Theta_\rho(\hat\theta)} J_\theta(\omega; x),
\qquad
J_\theta(\omega; x)
:=
\mathbb E_{z\sim\nu}\!\left[f\bigl(\omega,G_\theta(z,x)\bigr)\right],
\label{eq:robust_population_objective}
\end{equation}
and its empirical counterpart
\begin{equation}
\hat U_N(\omega;x)
:=
\sup_{\theta\in\Theta_\rho(\hat\theta)}\hat J_N(\omega;\theta,x),
\qquad
\hat J_N(\omega;\theta,x)
:=
\frac{1}{N}\sum_{i=1}^N f\bigl(\omega,G_\theta(z_i,x)\bigr).
\label{eq:sro_empirical_objective}
\end{equation}
The corresponding \emph{Sampler-Robust Optimization} (SRO) decision is
\begin{equation}
\hat{\omega}(x) \in \arg \min_{\omega \in \mathcal W} \hat U_N(\omega;x).
\label{eq:sro_decision}
\end{equation}

This formulation is directly aligned with generative decision pipelines: the optimizer evaluates decisions using simulated scenarios, and robustness is imposed on the very object that generates those scenarios.

\paragraph{Sharpness-aware interpretation.}
SRO also admits a natural sharpness-aware interpretation. Define the empirical and population $\theta$-sharpness by
\begin{align}
\hat{S}_\rho(\omega; x)
&= \sup_{\theta \in \Theta_\rho(\hat{\theta})}
\hat{J}_N(\omega; \theta, x)
- \hat{J}_N(\omega; x), \label{eq:empirical_sharpness} \\
S_\rho(\omega; x)
&= \sup_{\theta \in \Theta_\rho(\hat{\theta})}
J_\theta(\omega; x)
- J(\omega; x). \label{eq:population_sharpness}
\end{align}
Then the SRO rule can be written equivalently as
\begin{equation}
  \hat\omega(x)
  \;\in\;
  \arg\min_{\omega\in\mathcal W}
  \Bigl[
    \hat J_N(\omega;x)
    \;+\;
    \hat S_\rho(\omega;x)
  \Bigr].
  \label{eq:sharpness_form}
\end{equation}

Thus, SRO augments the nominal objective with a \emph{sampler-sharpness penalty}. Decisions that perform well only at the nominal sampler but deteriorate under nearby perturbations are penalized, whereas decisions whose performance remains stable across nearby samplers are favored. This is structurally reminiscent of sharpness-aware learning \citep{foret2021sharpness}, but the role of perturbation is different: in sharpness-aware learning, perturbations regularize model training; in SRO, they stress-test downstream decisions against sampler misspecification.

\paragraph{Probabilistic guarantee.}
The key insight is that robustification creates a \emph{population slack} that can partially absorb finite-simulation error. Define
\begin{equation}
J^\star(\omega;x) := J_{\theta^\star}(\omega;x),
\qquad
\mu(\omega,\rho;x)
:= U(\omega;x)-J^\star(\omega;x),
\qquad
\bar\mu(\rho;x)
:= \inf_{\omega\in\mathcal W}\mu(\omega,\rho;x).
\label{eq:population_slack}
\end{equation}
We say that the coverage condition holds at context $x$ if there exists some
$\tilde\theta^\star \in \Theta_\rho(\hat\theta)$ such that
\[
P_{\tilde\theta^\star}(\cdot\mid x)=P_{\theta^\star}(\cdot\mid x).
\]
Under this condition, $\mu(\omega,\rho;x)\ge 0$ for all $\omega\in\mathcal W$, and hence $\bar\mu(\rho;x)\ge 0$. The quantity $\bar\mu(\rho;x)$ captures the uniform margin by which the robust population objective dominates the true population objective.

The full theorem statement and proof are deferred to Appendix~\ref{app:generalization}. Informally, the main guarantee can be summarized as follows.
\begin{theorem}[Informal]
\label{thm:informal_main}
Suppose the coverage condition holds at context $x$. Under standard regularity conditions, with high probability over the latent samples $\{z_i\}_{i=1}^N$, the following bound holds uniformly over $\omega \in \mathcal W$:
\begin{equation}
J^\star(\omega; x)
\;\le\;
\hat U_N(\omega; x)
\;+\; \varepsilon_N(\mathcal W)
\;-\; \bar{\mu}(\rho; x),
\label{eq:informal_main_bound}
\end{equation}
where $\varepsilon_N(\mathcal W)$ is the finite-simulation error controlled by the complexity of the decision class, and $\bar{\mu}(\rho; x)$ is the population slack induced by robustification.
\end{theorem}

Thus, SRO addresses sampler misspecification and finite-simulation error through a single robustification mechanism: under coverage, the worst-case objective protects against misspecification, while the induced slack reduces the simulation error entering the certificate from $\varepsilon_N(\mathcal W)$ to $\varepsilon_N(\mathcal W)-\bar\mu(\rho;x)$. A notable feature of the guarantee is that the finite-simulation term $\varepsilon_N(\mathcal W)$ depends on the complexity of the decision class $\mathcal W$, but not on the parameter dimension or covering complexity of the generator neighborhood $\Theta_\rho(\hat\theta)$. This keeps the certificate statistically meaningful even for large, highly parameterized generators. The proof, given in Appendix~\ref{app:generalization}, shows how this follows from the one-sided structure of the robust certificate.

\paragraph{Three regimes of the radius $\rho$.}
The effect of \eqref{eq:informal_main_bound} depends on the size of the induced slack $\bar\mu(\rho;x)$:
\begin{enumerate}
\item
\emph{Small $\rho$} $(\bar\mu(\rho;x)\approx 0)$:
the bound essentially reduces to
\[
J^\star(\omega;x)\le\hat U_N(\omega;x)+\varepsilon_N,
\]
so robustification offers little additional protection against finite-simulation error.

\item
\emph{Moderate $\rho$} $(0<\bar\mu(\rho;x)<\varepsilon_N)$:
the finite-simulation correction is reduced, as the conservatism introduced by robustification absorbs part of the Monte Carlo error.

\item
\emph{Large $\rho$} $(\bar\mu(\rho;x)\ge\varepsilon_N)$:
the correction is fully absorbed, so $\hat U_N(\omega;x)$ itself serves as a valid upper certificate for $J^\star(\omega;x)$.
\end{enumerate}

These regimes reveal a fundamental tradeoff governed by $\rho$: larger perturbation neighborhoods increase conservatism, but can also strengthen the resulting reliability certificate.

% \begin{remark}[When can $\bar\mu(\rho;x)>0$?]
% \label{rem:slack_positive}
% In most nondegenerate settings, one expects $\bar\mu(\rho;x)>0$.
% \textcolor{red}{Indeed, suppose there exists some $\tilde\theta^\star \in \Theta_\rho(\hat\theta)$ such that
% \[
% P_{\tilde\theta^\star}(\cdot\mid x)=P_{\theta^\star}(\cdot\mid x).
% \]
% If $\tilde\theta^\star$ lies in the interior of $\Theta_\rho(\hat\theta)$ and the map
% $\theta\mapsto J_\theta(\omega;x)$ varies nontrivially near
% $\tilde\theta^\star$ in a manner that is sufficiently uniform over
% $\omega\in\mathcal W$, robustification creates a strict separation
% between the worst-case objective and the true objective across the
% decision class.}

% By contrast, $\bar\mu(\rho;x)=0$ corresponds to a
% degenerate alignment in which this separation vanishes at the
% decision(s) that matter most---that is, robustification fails to
% improve upon the true objective in the least favorable direction.
% Thus, the key issue is not merely pointwise variation in
% $\theta\mapsto J_\theta(\omega;x)$, but whether such variation yields
% a nontrivial margin uniformly over $\omega\in\mathcal W$.
% \end{remark}

%%%%%%%%%%%%%%%%%%%%%%%%%%%%%%%%%%%%%%%%%%%%%%
\section{Solving the Sampler-Robust Problem}
\label{sec:solver}

The sampler-robust formulation in \eqref{eq:sro_decision} leads to a constrained adversarial problem over both the decision variable and the sampling mechanism. In general, the inner maximization over sampler parameters is nonlinear and may be high-dimensional, so exact solution can be difficult. We therefore develop practical optimization procedures tailored to two natural computational regimes. When the uncertainty radius $\rho$ is small, the inner adversary can be well approximated by a first-order local perturbation around the nominal sampler. When $\rho$ is moderate or large, such a local approximation may become inaccurate, and it is preferable to optimize the constrained min--max problem more directly through alternating projected updates.

\subsection{First-order worst-case sampler refinement}

We first consider the regime in which the ambiguity radius $\rho$ is small, so that the inner maximization can be approximated locally around the nominal generator parameter $\hat\theta$. Let
\[
\varepsilon = \theta - \hat{\theta}.
\]
Then, for a fixed decision $\omega$, the inner problem in \eqref{eq:sro_decision} takes the form
\begin{equation}
\max_{\|\varepsilon\|_p \le \rho}
\frac{1}{N}\sum_{i=1}^N f\bigl(\omega, G_{\hat{\theta}+\varepsilon}(z_i,x)\bigr),
\label{eq:first_order_inner}
\end{equation}
where $\|\cdot\|_p$ is a chosen norm on parameter space.

When $\rho$ is small, a first-order Taylor expansion around $\hat\theta$ yields
\begin{equation}
\frac{1}{N}\sum_{i=1}^N f\bigl(\omega, G_{\hat{\theta}+\varepsilon}(z_i,x)\bigr)
\approx
\hat{J}_N(\omega;x) + \varepsilon^\top g,
\label{eq:first_order_taylor}
\end{equation}
where
\begin{equation}
g
:=
\nabla_\theta
\left(
\frac{1}{N}\sum_{i=1}^N f\bigl(\omega, G_{\theta}(z_i,x)\bigr)
\right)\Bigg|_{\theta=\hat{\theta}}.
\label{eq:first_order_gradient}
\end{equation}
The approximate inner problem therefore reduces to the linear maximization
\begin{equation}
\varepsilon^\star
\in
\arg\max_{\|\varepsilon\|_p \le \rho}
\varepsilon^\top g.
\label{eq:dual_inner_problem}
\end{equation}
For $1<p<\infty$ and $g\neq 0$, dual norm theory gives
\[
\max_{\|\varepsilon\|_p\le \rho} \varepsilon^\top g
=
\rho\|g\|_q,
\qquad \frac1p+\frac1q=1,
\]
with optimizer
\[
\varepsilon^\star
=
\rho\,
\frac{\operatorname{sign}(g)\odot |g|^{q-1}}
{\|g\|_q^{q-1}}.
\]
The resulting perturbed parameter
\[
\tilde{\theta} = \hat{\theta} + \varepsilon^\star
\]
can be interpreted as a locally worst-case sampler for the current decision iterate.

% For implementation, given latent draws $\{z^{(i)}\}_{i=1}^B \sim \nu$, define
% \begin{equation}
% J_B(\omega,\theta)
% :=
% \frac{1}{B}\sum_{i=1}^{B} f\bigl(\omega, G_{\theta}(z^{(i)},x)\bigr).
% \label{eq:minibatch_objective}
% \end{equation}

Algorithm~\ref{alg:first_order_sampler} uses this local adversarial approximation to refine the decision variable. At each iteration, it computes\footnote{In our implementation, this latent batch is sampled once and then reused throughout the refinement procedure. This keeps the simulation budget fixed and ensures a fair comparison with the nominal sample-average baseline, which is optimized using the same number of generated scenarios.} the worst-case first-order perturbation around $\hat\theta$ for the current decision iterate, then updates the decision against the resulting perturbed sampler.

\begin{algorithm}[t]
\caption{First-Order Worst-Case Sampler Refinement}
\label{alg:first_order_sampler}
\small
\begin{algorithmic}[1]
\Require Nominal generator $G_{\hat{\theta}}$, context $x$, feasible set $\mathcal{W}$, batch size $B$, radius $\rho$, norm $\|\cdot\|_p$ with dual exponent $q$, refinement steps $K$, step size $\alpha_{\omega}$
\Ensure Robust decision $\omega^{\mathrm{rob}}$

\Statex \textit{Mini-batch objective:} For latent draws $\{z^{(i)}\}_{i=1}^{B}$, define
\[
J_B(\omega,\theta)
\triangleq
\frac{1}{B}\sum_{i=1}^{B} f\!\left(\omega, G_{\theta}(z^{(i)}, x)\right).
\]

\State Sample i.i.d.\ latent noise $\{z^{(i)}\}_{i=1}^{B} \sim \nu$
\State Initialize $\omega^{(0)} \in \mathcal{W}$

\For{$k = 0$ to $K-1$}
    \State $g^{(k)} \gets \nabla_{\theta} J_B(\omega^{(k)},\theta)\big|_{\theta=\hat{\theta}}$
    \If{$\|g^{(k)}\|_q = 0$}
        \State $\varepsilon^{\star,(k)} \gets 0$
    \Else
        \State $\varepsilon^{\star,(k)}
        \gets
        \rho\,
        \dfrac{\mathrm{sign}(g^{(k)})\odot |g^{(k)}|^{q-1}}
        {\|g^{(k)}\|_q^{\,q-1}}$
    \EndIf
    \State $\tilde{\theta}^{(k)} \gets \hat{\theta} + \varepsilon^{\star,(k)}$
    \State $\omega^{(k+1)}
    \gets
    \Pi_{\mathcal{W}}\!\left(
    \omega^{(k)} - \alpha_{\omega}
    \nabla_{\omega} J_B(\omega,\tilde{\theta}^{(k)})\big|_{\omega=\omega^{(k)}}
    \right)$
\EndFor

\State \Return $\omega^{\mathrm{rob}} \gets \omega^{(K)}$
\end{algorithmic}
\end{algorithm}

%This procedure is computationally attractive because the adversarial refinement is obtained in closed form once the gradient with respect to $\theta$ is available. Its main limitation is that it relies on local geometry around $\hat\theta$, and is therefore best suited to settings in which the worst-case sampler lies near the nominal one.

\subsection{Projected alternating minimax refinement}

When the uncertainty radius $\rho$ is moderate or large, the local linearization in \eqref{eq:first_order_taylor} may no longer capture the inner adversary adequately. In that regime, we instead optimize the sampler-robust objective more directly via alternating projected descent--ascent steps applied to the fixed-batch objective
\begin{equation}
\min_{\omega \in \mathcal{W}}
\max_{\theta \in \Theta_{\rho}(\hat{\theta})}
J_B(\omega,\theta).
\label{eq:alternating_minimax}
\end{equation}

We adopt a projected gradient descent--ascent scheme that alternates:
(i) a projected ascent step in $\theta$ to refine the adversarial sampler within $\Theta_\rho(\hat{\theta})$, and
(ii) a projected descent step in $\omega$ to adapt the decision variable to the updated worst-case sampler. Because the generator-parameter maximization is typically nonconcave, aggressive adversarial updates can destabilize the minimax dynamics. We therefore adopt a stabilized two-timescale scheme with $\alpha_\theta \ll \alpha_\omega$: the adversarial sampler is updated cautiously, while the decision variable, often governed by a convex or better-conditioned problem, is allowed to adapt more rapidly. This choice follows the broader two-timescale principle that asymmetric update speeds can stabilize minimax dynamics, while reflecting the particular asymmetry of our setting. Algorithm~\ref{alg:two_timescale_sampler} summarizes the procedure.

\begin{algorithm}[t]
\caption{Two-Timescale Alternating Minimax Refinement}
\label{alg:two_timescale_sampler}
\small
\begin{algorithmic}[1]
\Require Nominal generator $G_{\hat{\theta}}$, context $x$, feasible set $\mathcal{W}$, batch size $B$, radius $\rho$, norm $\|\cdot\|_p$, refinement steps $K$, step sizes $\alpha_{\theta}, \alpha_{\omega}$ with $\alpha_{\theta} \ll \alpha_{\omega}$
\Ensure Robust decision $\omega^{\mathrm{rob}}$

\Statex \textit{Mini-batch objective:} For latent draws $\{z^{(i)}\}_{i=1}^{B}$, define
\[
J_B(\omega,\theta)
\triangleq
\frac{1}{B}\sum_{i=1}^{B} f\!\left(\omega, G_{\theta}(z^{(i)}, x)\right).
\]

\State Sample i.i.d.\ latent noise $\{z^{(i)}\}_{i=1}^{B} \sim \nu$
\State Initialize $\theta^{(0)} \gets \hat{\theta}$ and choose any $\omega^{(0)} \in \mathcal{W}$

\For{$k = 0$ to $K-1$}
    \State \textit{Adversary step (slow timescale):}
    \State $\theta' \gets \theta^{(k)} + \alpha_{\theta}\nabla_{\theta} J_B(\omega^{(k)},\theta)\big|_{\theta=\theta^{(k)}}$
    \State $\theta^{(k+1)} \gets \hat{\theta} + \Pi_{\|\cdot\|_p \le \rho}\!\left(\theta' - \hat{\theta}\right)$

    \State \textit{Decision step (fast timescale):}
    \State $\omega^{(k+1)} \gets \Pi_{\mathcal{W}}\!\left(
    \omega^{(k)} - \alpha_{\omega}\nabla_{\omega} J_B(\omega,\theta^{(k+1)})\big|_{\omega=\omega^{(k)}}
    \right)$
\EndFor

\State \Return $\omega^{\mathrm{rob}} \gets \omega^{(K)}$
\end{algorithmic}
\end{algorithm}

\paragraph{Discussion.}
The two procedures above serve complementary roles. Algorithm~\ref{alg:first_order_sampler} provides an efficient local refinement when the worst-case sampler is well approximated by first-order geometry around the nominal parameter. Algorithm~\ref{alg:two_timescale_sampler} offers a more direct exploration of the constrained adversarial problem when the ambiguity set is larger or the local approximation is less reliable. Together, they provide practical optimization tools for sampler-robust optimization across a range of operating regimes.

%%%%%%%%%%%%%%%%%%%%%%%%%%%%%%%%%%%%%%%%%%%%%%%%%%%%

%%%%%%%%%%%%%%%%%%%
%%%%%%%%%%%%%%%%%%
\section{Numerical Experiments}
\label{sec:experiments}

We study the proposed sampler-robust framework in a portfolio optimization setting. The experiments serve two purposes. The first is to test the reliability mechanism suggested by the Theorem~\ref{thm:informal_main} in a controlled setting where both the oracle sampler and the nominal sampler are observable. The second is to examine whether the same robustness effect translates into improved out-of-sample portfolio performance on real market data.

\subsection{Portfolio optimization: formulation and experimental setup}

At each decision time $t$, the conditioning state is a rolling window of the past
$L=10$ daily multivariate log returns,
\[
x_t \in \mathbb{R}^{L \times d},
\]
where $d$ is the number of assets in the portfolio universe. For a generator
parameter $\theta$, the conditional generator $G_\theta(\cdot, x_t)$ maps a latent
noise vector together with the historical return window $x_t$ to a simulated
one-step-ahead return vector in $\mathbb{R}^d$. A portfolio is a weight vector
\[
\omega \in \mathcal{W}
:=
\left\{
\omega\in\mathbb{R}^d:
\sum_{j=1}^d \omega_j = 1,\;
\omega_j \ge 0
\right\},
\]
so that $\mathcal W$ is the long-only simplex.

For a fixed state $x_t$ and generator parameter $\theta$, let
\[
\tilde r_\theta^{(1)},\dots,\tilde r_\theta^{(N)}
\]
denote Monte Carlo draws from the conditional generator $G_\theta(\cdot,x_t)$.
The portfolio return in scenario $i$ is
\begin{equation}
\pi^{(i)}(\omega,\theta;x_t)
:=
(\tilde r_\theta^{(i)})^\top \omega.
\label{eq:generated_portfolio_return}
\end{equation}

Throughout the experiments, we use the quadratic utility
\begin{equation}
u(\pi)=\pi-\frac{\lambda}{2}\pi^2,
\label{eq:quadratic_utility}
\end{equation}
where $\lambda>0$ is the risk-aversion parameter. We fix $\lambda=10$ in all experiments. The empirical utility is
\begin{equation}
\hat J_N(\omega;\theta,x_t)
=
\frac{1}{N}\sum_{i=1}^N
\left[
\pi^{(i)}(\omega,\theta;x_t)
-\frac{\lambda}{2}\bigl(\pi^{(i)}(\omega,\theta;x_t)\bigr)^2
\right].
\label{eq:general_empirical_utility}
\end{equation}

Because the experiments are written under a utility convention rather than the loss convention used in the general formulation, the optimization direction is reversed: the decision maker maximizes utility and the adversarial sampler minimizes it.
The nominal portfolio is obtained by solving \eqref{eq:general_empirical_utility}
under the fitted generator parameter $\hat\theta$. The sampler-robust portfolio
instead optimizes the corresponding worst-case empirical utility over the parameter ball $\Theta_\rho(\hat\theta)$.
% \begin{equation}
% \Theta_\rho(\hat\theta)
% :=
% \{\theta:\|\theta-\hat\theta\|\le\rho\},
% \label{eq:robust_empirical_portfolio_obj}
% \end{equation}
% with robustness radius $\rho>0$.

In both experiments, we use the $\ell_2$ norm to define the parameter neighborhood, corresponding to the Euclidean perturbation geometry commonly used in sharpness-aware optimization. In addition to the training data, we reserve 200 held-out observations for model selection and evaluation. The first 100 held-out observations form a validation window used to select the robustness radius, while the remaining 100 are used for out-of-sample testing. The reported results use $\rho=0.3$, selected by this validation procedure and then fixed throughout the test-period comparison.

\subsection{Nominal generative model}

As the nominal generative model, we use a conditional LSTM-GAN trained on rolling windows of multivariate log returns. The architecture is inspired by Fin-GAN-style generative models for financial time series \citep{vuletic2024fin}. The generator combines a latent noise vector with the historical return window and outputs a one-step-ahead multivariate return vector, while the discriminator takes the conditioning sequence together with either the realized or generated next-period return and produces a scalar score. Throughout the experiments, we fix the latent noise dimension at $d_z=8$ and the LSTM hidden dimension at $D_h=8$. Training is based on a least-squares GAN objective. Detailed architectural specifications are deferred to Appendix~\ref{subsec:architecture}. The same nominal architecture is used throughout so that performance differences can be attributed to sampler-robust optimization rather than to changes in the predictive model.

\paragraph{Generator validity screening.}
Because GAN training can be unstable and susceptible to mode collapse, we impose a post-training validity screen before conducting downstream portfolio evaluation in both experiments \citep{arjovsky2017wasserstein}. This screen is based on basic distributional diagnostics comparing generated and observed return samples. Runs exhibiting clear evidence of collapse or other degeneracies are excluded, since in such cases the fitted generator does not provide a reliable basis for economically meaningful portfolio evaluation. The reported results therefore compare nominal and sampler-robust portfolio construction after excluding failed generator-training runs.

\subsection{Diagnostic and evaluation metrics}\label{metrics}

To illustrate the reliability mechanism suggested by Theorem~\ref{thm:informal_main}, we introduce several diagnostic quantities for the controlled experiment. These diagnostics evaluate whether the optimized empirical criterion remains close to oracle performance under the oracle generator. In addition, to assess realized portfolio performance, we report standard out-of-sample return and risk metrics.

\paragraph{Controlled-experiment diagnostics.}
In the controlled experiment, where the oracle generator is known, we report empirical utility, oracle utility, and the empirical-to-oracle gap. For any portfolio $\omega$, the empirical utility is the optimization criterion evaluated under the generator used by the decision rule. For the nominal method, this is
\[
\hat J_N(\omega;\hat\theta,x_t),
\]
while for the robust method it is the robust empirical objective
\[
\hat U_N(\omega;x_t).
\]

Let $\theta^\star$ denote the oracle generator parameter. The oracle utility of a portfolio $\omega$ is defined as
\begin{equation}
J^\star(\omega;x_t)
=
\mathbb E_{z\sim\nu}
\left[
u\!\left(G_{\theta^\star}(z,x_t)^\top \omega\right)
\right].
\label{eq:oracle_utility_def}
\end{equation}
In practice, this quantity is approximated using a large Monte Carlo sample from the oracle generator.

For the nominal portfolio, the empirical-to-oracle gap is defined by
\begin{equation}
\mathrm{Gap}^{\mathrm{nom}}_t
:=
\hat J_N(\omega_t^{\mathrm{nom}};\hat\theta,x_t)
-
J^\star(\omega_t^{\mathrm{nom}};x_t),
\label{eq:nominal_gap_def}
\end{equation}
and for the robust portfolio by
\begin{equation}
\mathrm{Gap}^{\mathrm{rob}}_t
:=
\hat U_N(\omega_t^{\mathrm{rob}};x_t)
-
J^\star(\omega_t^{\mathrm{rob}};x_t).
\label{eq:robust_gap_def}
\end{equation}
These gap measures quantify how closely the optimized empirical criterion reflects oracle performance.

\paragraph{Out-of-sample performance metrics.}
For out-of-sample evaluation, let $r_{t+1}\in\mathbb{R}^d$ denote the realized next-period log-return vector. The realized net portfolio return is
\begin{equation}
R_t(\omega_t)
=
\omega_t^\top \exp(r_{t+1}) - 1,
\label{eq:oos_measurement}
\end{equation}
where the exponential is taken componentwise. Given the realized return sequence $\{R_t\}_{t=1}^T$, we report mean return, standard deviation, Sharpe ratio, CVaR(5\%), and maximum drawdown.
\subsection{Experimental design}

In both experiments, we compare the proposed sampler-robust portfolio with its nominal counterpart, which solves the same portfolio problem under the fitted generator without robust refinement. Since the two methods share the same conditional LSTM-GAN, scenario budget, and utility specification, the comparison isolates the contribution of sampler-level robustness. In the controlled experiment, we also report an oracle portfolio computed under the oracle generator. This portfolio is not implementable in practice and is included only as a diagnostic benchmark.

Both experiments are repeated over 20 random seeds, indexed by 40, 41, \dots, 59. In the controlled experiment, all seeded runs are retained and no evidence of generator collapse is detected. In the real-data experiment, only seed 49 fails the post-training validity screen and is excluded from the downstream portfolio comparison. This distinction is consistent with the greater difficulty of training conditional generators on real return data than in the controlled synthetic setting.

\subsubsection{Controlled generator-to-generator experiment}

The first experiment studies the reliability mechanism suggested by Theorem~\ref{thm:informal_main} in a controlled setting where both the oracle and nominal samplers are observable. We fix the 20-stock universe
\[
\left\{
\begin{aligned}
&\text{AAPL}, \text{MSFT}, \text{GOOGL}, \text{AMZN}, \text{META}, \\
&\text{TSLA}, \text{NVDA}, \text{ADBE}, \text{INTC}, \text{CRM}, \\
&\text{AMD}, \text{CSCO}, \text{ORCL}, \text{IBM}, \text{QCOM}, \\
&\text{JPM}, \text{BAC}, \text{WFC}, \text{C}, \text{GS}
\end{aligned}
\right\}.
\]
A generator trained on the historical data is treated as the oracle sampler. We then simulate a synthetic future path of length 1000, retrain a nominal generator on the initial 800 observations, calibrate the robustness radius on the subsequent 100 observations, and assess portfolio decisions over the final 100-observation test period. At each decision time in the test period, both the nominal and sampler-robust portfolios are computed using the retrained generator. Because the oracle generator remains observable in this controlled setting, we can additionally evaluate oracle utility and the empirical-to-oracle gap defined in Section~\ref{metrics}.

\subsubsection{Real-data experiment with random stock pools}

The second experiment evaluates practical robustness on real market data. In each seeded run, we
randomly select a 10-stock universe from a panel of 50 S\&P 500 stocks, with additional details on the
stock universe provided in the Appendix \ref{subsec:stock_universe}, and train a conditional LSTM-GAN using rolling windows constructed from the first 1000 observations. We then reserve the
subsequent 200 observations, using the first 100 as a validation window for selecting the robustness
radius and the final 100 as the out-of-sample test period. Portfolio decisions are made sequentially over
this 100-day test window.
%This design is used to maintain stable nominal generator training across seeds. In preliminary experiments,larger stock universes led to less reliable GAN training and a higher incidence of collapse, which undermines the interpretability of the downstream portfolio comparison.
This random-pool design also reduces dependence
on any particular hand-picked universe and tests whether the method remains effective across different
cross-sectional compositions and training noise. In this setting, no oracle generator is available, so evaluation
is based on realized out-of-sample performance alone. 

%Whereas the controlled study evaluates the reliability mechanism implied by the theory, the real-data study evaluates its practical portfolio implications.

\subsection{Results}
\label{subsec:results}

\subsubsection{Controlled generator-to-generator results}

\begin{figure}[H]
    \centering
    \includegraphics[width=0.92\textwidth]{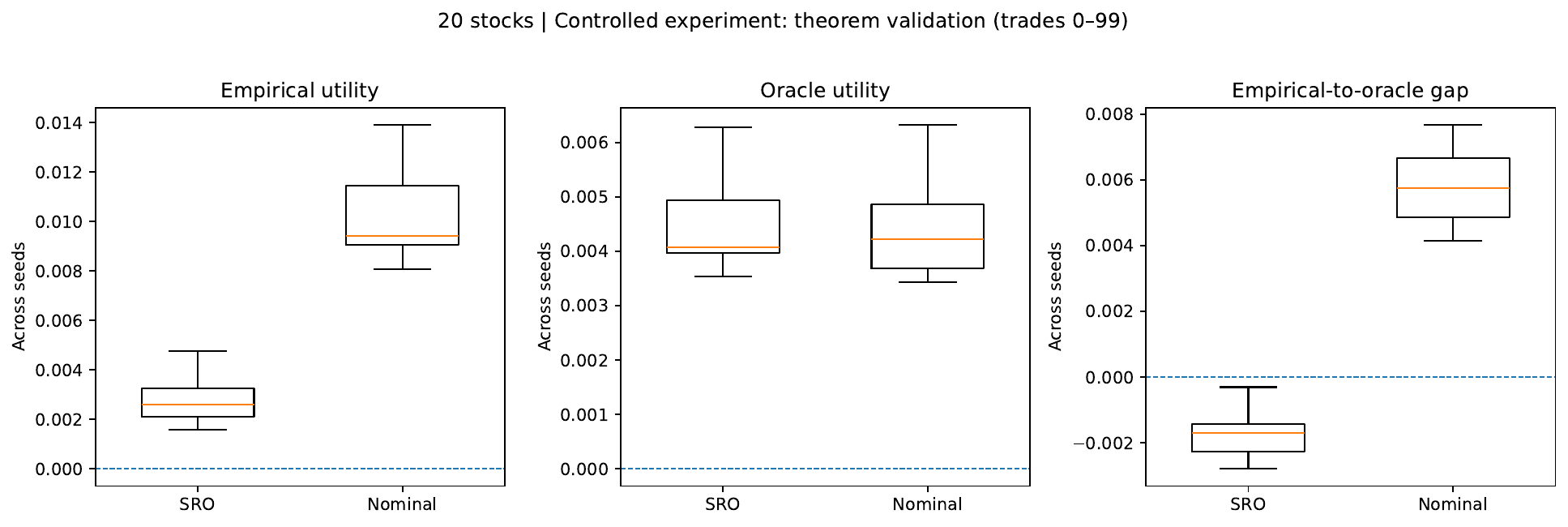}
    \caption{Utility-based diagnostics in the controlled experiment: empirical utility, oracle utility, and the empirical-to-oracle gap.}
    \label{fig:controlled_theorem_validation}
\end{figure}

\begin{figure}[H]
    \centering
    \includegraphics[width=0.88\textwidth]{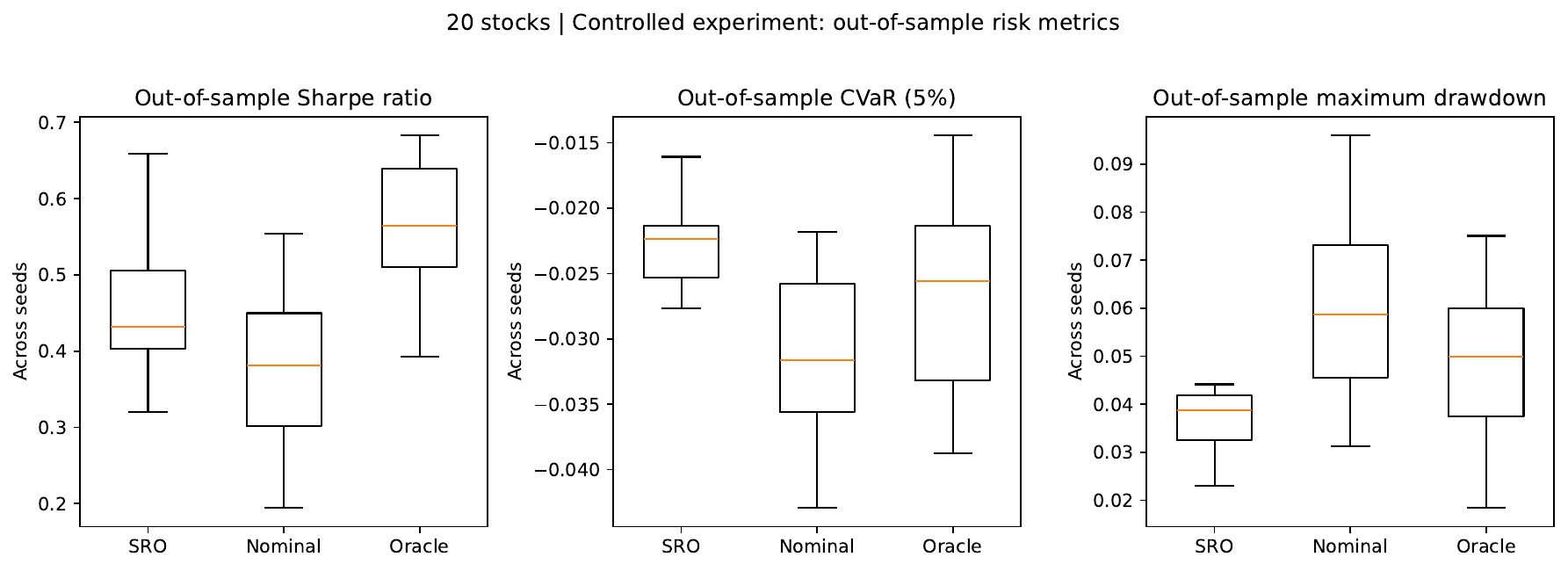}
    \caption{Out-of-sample risk metrics in the controlled experiment: Sharpe ratio, CVaR(5\%), and maximum drawdown.}
    \label{fig:controlled_oos_risk}
\end{figure}

Figure~\ref{fig:controlled_theorem_validation} presents the utility-based diagnostics in the controlled experiment, including empirical utility, oracle utility, and the empirical-to-oracle gap. Figure~\ref{fig:controlled_oos_risk} reports the corresponding out-of-sample risk metrics, while Table~\ref{tab:controlled_summary} summarizes the main numerical results across 20 seeds.

\begin{table}[H]
\centering
\caption{Controlled generator-to-generator experiment: summary across 20 seeds}
\label{tab:controlled_summary}
\small
\setlength{\tabcolsep}{6pt}
\renewcommand{\arraystretch}{1.15}

\begin{tabular}{lccc}
\toprule
& \multicolumn{3}{c}{\textbf{(a) Utility alignment}} \\
\cmidrule(lr){2-4}
\textbf{Method} & \textbf{Empirical Utility} & \textbf{Oracle Utility} & \textbf{Empirical-to-Oracle Gap} \\
\midrule
SRO
& $0.002814 \pm 0.001013$
& $0.004483 \pm 0.000821$
& $\mathbf{-0.001669 \pm 0.000730}$ \\

Nominal
& $\mathbf{0.010258 \pm 0.001837}$
& $\mathbf{0.004493 \pm 0.001021}$
& $0.005765 \pm 0.001092$ \\
\midrule
\multicolumn{4}{c}{} \\
\multicolumn{4}{c}{\textbf{(b) Out-of-sample risk}} \\
\cmidrule(lr){2-4}
\textbf{Method} & \textbf{Sharpe Ratio} & \textbf{CVaR(5\%)} & \textbf{Maximum Drawdown} \\
\midrule
SRO
& $0.455701 \pm 0.094245$
& $\mathbf{-0.022958 \pm 0.004727}$
& $\mathbf{0.038815 \pm 0.011004}$ \\

Nominal
& $0.374830 \pm 0.101227$
& $-0.032756 \pm 0.008279$
& $0.063319 \pm 0.023607$ \\

Oracle
& $\mathbf{0.578916 \pm 0.103024}$
& $-0.026248 \pm 0.006919$
& $0.051941 \pm 0.021671$ \\
\bottomrule
\end{tabular}
\end{table}

Several patterns emerge clearly. First, the nominal portfolio attains the highest empirical utility. This is expected, since it is optimized directly under the fitted nominal sampler and therefore reflects the most optimistic empirical criterion. By contrast, the robust portfolio is more conservative at the empirical level.

Second, this empirical advantage translates only weakly into oracle performance. The oracle utilities of the robust and nominal portfolios remain very close, while the empirical-to-oracle gap is dramatically smaller under SRO. In particular, the average gap is close to zero and slightly negative for SRO, but clearly positive for the nominal portfolio. Thus, the robust criterion is markedly less optimistic relative to oracle performance. This pattern is consistent with the mechanism suggested by the theory: robustification creates slack between the empirical criterion and the oracle utility, thereby reducing effective over-optimism.

Third, the reduction in optimism is accompanied by improved out-of-sample risk control. Relative to the nominal portfolio, SRO achieves a higher Sharpe ratio, a less adverse CVaR(5\%), and a smaller maximum drawdown. The oracle portfolio attains the highest Sharpe ratio overall, as expected for a diagnostic benchmark constructed under the oracle generator, but SRO delivers the strongest drawdown control among the implementable methods. Hence, in the controlled setting, tighter alignment between the optimization criterion and oracle performance is associated with economically meaningful gains in downside stability.
\subsubsection{Real-data backtest results}
Tables~\ref{tab:real_return_riskadjusted}--\ref{tab:real_tailrisk} report the real-data backtest results across valid random seeds. 

%In the original set of runs, seed 49 was excluded because the trained generator failed a post-training quality check and exhibited collapse, producing degenerate simulated return scenarios. Since both the nominal and robust portfolio decisions are defined with respect to the learned generator, such a failed run does not yield a meaningful downstream comparison. The summary statistics below are therefore computed over the remaining valid seeds.

\begin{table}[H]
\centering
\caption{Real-data backtest: return and risk-adjusted performance}
\label{tab:real_return_riskadjusted}
\small
\setlength{\tabcolsep}{8pt}
\renewcommand{\arraystretch}{1.15}

\begin{tabular}{lccc}
\toprule
\textbf{Method} & \textbf{Mean Return} & \textbf{Std. Dev.} & \textbf{Sharpe Ratio} \\
\midrule
SRO
& $\mathbf{0.000565 \pm 0.000802}$
& $\mathbf{0.013447 \pm 0.004547}$
& $\mathbf{0.046313 \pm 0.064395}$ \\

Nominal
& $0.000508 \pm 0.000955$
& $0.018575 \pm 0.005747$
& $0.032855 \pm 0.057811$ \\
\bottomrule
\end{tabular}
\end{table}

\begin{table}[H]
\centering
\caption{Real-data backtest: downside and drawdown risk}
\label{tab:real_tailrisk}
\small
\setlength{\tabcolsep}{8pt}
\renewcommand{\arraystretch}{1.15}

\begin{tabular}{lcc}
\toprule
\textbf{Method} & \textbf{CVaR(5\%)} & \textbf{Maximum Drawdown} \\
\midrule
SRO
& $\mathbf{-0.029920 \pm 0.011380}$
& $\mathbf{0.101745 \pm 0.049559}$ \\

Nominal
& $-0.042013 \pm 0.013278$
& $0.148037 \pm 0.066574$ \\
\bottomrule
\end{tabular}
\end{table}
The real-data results across valid seeds show that SRO delivers more stable portfolio performance than the nominal benchmark. In particular, SRO achieves a higher mean return, lower volatility, a higher Sharpe ratio, and a smaller maximum drawdown. It also exhibits improved downside-risk behavior overall.

Taken together, these findings suggest that, conditional on a valid trained generator, the robust portfolio construction provides a more reliable risk-return tradeoff than the nominal alternative in the real-data setting.

%%%%%%%%%%%%%%%%%%%%%%%%%%%%%%%%%%%%%%%%%%%%%%%%%%%%
\section{Conclusion}

This paper proposed \emph{Sampler-Robust Optimization} (SRO), a robust decision-making framework for stochastic optimization problems driven by learned generative models. We showed that sampler-level robustification admits a sharpness-aware interpretation and yields a high-probability reliability certificate under distributional coverage, explaining why robustifying the sampler can improve decision performance even after the nominal generator has been trained. Portfolio-optimization experiments show that SRO produces more stable decisions and stronger out-of-sample performance under distribution shift. Overall, our results highlight the value of treating the sampler as the operational object of uncertainty in simulation-based decision pipelines.

\appendix

\section{Formal Reliability Guarantee and Proof}
\label{app:generalization}

Here, we state and prove the formal probabilistic guarantee underlying the informal bound in Section~\ref{sec:formulation}. For clarity, we present the result for Gaussian latent noise, which allows us to invoke standard concentration inequalities for Lipschitz functions of Gaussian random variables. The same argument extends to other latent laws that admit analogous concentration bounds.

\subsection{Assumptions and setup}

Fix a context $x\in\mathcal X$. Throughout this appendix, we work under the generative representation
\[
y = G_\theta(z,x),
\qquad
z\sim\mathcal N(0,I),
\]
and define, for $\omega\in\mathcal W$ and $\theta\in\Theta$,
\[
\ell_\theta(z,\omega)
:=
f\bigl(\omega,G_\theta(z,x)\bigr).
\]
Recall that the true population objective is
\[
J^\star(\omega;x)
=
\mathbb E_{z\sim\mathcal N(0,I)}[\ell_{\theta^\star}(z,\omega)],
\]
while the empirical sampler-robust objective is
\[
\hat U_N(\omega;x)
=
\sup_{\theta\in\Theta_\rho(\hat\theta)}
\frac{1}{N}\sum_{i=1}^N \ell_\theta(z_i,\omega),
\qquad
z_1,\dots,z_N \stackrel{\mathrm{iid}}{\sim}\mathcal N(0,I).
\]
Likewise,
\[
U(\omega;x)
=
\sup_{\theta\in\Theta_\rho(\hat\theta)}
J_\theta(\omega;x),
\qquad
J_\theta(\omega;x)
=
\mathbb E_{z\sim\mathcal N(0,I)}[\ell_\theta(z,\omega)].
\]

For each $\omega\in\mathcal W$, define the pointwise slack
\[
\mu(\omega,\rho;x)
:=
U(\omega;x)-J^\star(\omega;x),
\]
and the population slack
\[
\bar\mu(\rho;x)
:=
\inf_{\omega\in\mathcal W}\mu(\omega,\rho;x).
\]

We impose the following assumptions:
\begin{itemize}
    \item[(A0)] \textbf{Distributional coverage.} There exists some
    \[
    \tilde\theta^\star \in \Theta_\rho(\hat\theta)
    \]
    such that
    \[
    P_{\tilde\theta^\star}(\cdot\mid x)=P_{\theta^\star}(\cdot\mid x).
    \]

    \item[(A1)] \textbf{Metric entropy of the decision set.} For every $\varepsilon>0$, the decision set $\mathcal W$ admits an $\varepsilon$-cover of size $N_{\mathcal W}(\varepsilon)$.

    \item[(A2)] \textbf{Lipschitz continuity in the decision variable.} There exists $L_\omega>0$ such that for all $y\in\mathcal Y$ and all $\omega,\omega'\in\mathcal W$,
    \[
    |f(\omega,y)-f(\omega',y)|
    \le
    L_\omega\|\omega-\omega'\|.
    \]

    \item[(A3)] \textbf{Lipschitz continuity in the outcome variable.} There exists $L_y>0$ such that for all $\omega\in\mathcal W$ and all $y,y'\in\mathcal Y$,
    \[
    |f(\omega,y)-f(\omega,y')|
    \le
    L_y\|y-y'\|.
    \]

    \item[(A4)] \textbf{Uniform Lipschitz continuity of the generator in the latent variable.} There exists $L_z>0$ such that for all $\theta\in\Theta_\rho(\hat\theta)$,
    \[
    \|G_{\theta}(z,x)-G_{\theta}(z',x)\|
    \le
    L_z\|z-z'\|
    \qquad
    \text{for all } z,z'\in\mathcal Z.
    \]

    \item[(A5)] \textbf{Attainment of the worst-case population sampler.} For every $\omega\in\mathcal W$, the supremum defining $U(\omega;x)$ is attained; that is, there exists
    \[
    \bar\theta(\omega)\in\Theta_\rho(\hat\theta)
    \quad\text{such that}\quad
    U(\omega;x)=J_{\bar\theta(\omega)}(\omega;x).
    \]
\end{itemize}

Under (A3) and (A4), the map
\[
z\mapsto f(\omega,G_\theta(z,x))
\]
is $L_yL_z$-Lipschitz uniformly over $\omega\in\mathcal W$ and $\theta\in\Theta_\rho(\hat\theta)$. This is the key regularity needed for concentration of the empirical objective around its population counterpart.

\subsection{Formal theorem}

\begin{theorem}[Uniform bound under coverage]
\label{thm:sro_generalization}
Suppose (A0)--(A5) hold. Then for any $\varepsilon>0$ and any $\delta\in(0,1)$, conditional on $\hat\theta$, $\theta^\star$, and the context $x$, with probability at least $1-\delta$ over the latent draws $z_1,\dots,z_N \stackrel{\mathrm{iid}}{\sim}\mathcal N(0,I)$, the following bound holds simultaneously for all $\omega\in\mathcal W$:
\begin{equation}
J^\star(\omega;x)
\le
\hat U_N(\omega;x)
+
L_yL_z
\sqrt{
\frac{
2\bigl(\log N_{\mathcal W}(\varepsilon)+\log(1/\delta)\bigr)
}{N}}
-
\bar\mu(\rho;x)
+
2L_\omega\varepsilon.
\label{eq:sro_generalization_bound}
\end{equation}
\end{theorem}

The certificate in \eqref{eq:sro_generalization_bound} has the same overall structure as a standard uniform concentration bound, but with two important refinements. First, robustification produces the slack term $-\bar\mu(\rho;x)$, where $\bar\mu(\rho;x)\ge 0$ under distributional coverage, so the same mechanism that guards against sampler misspecification also partially offsets the finite-simulation correction. Second, the finite-simulation term depends on the complexity of the decision class $\mathcal W$, but not on the parameter dimension or covering complexity of the generator neighborhood $\Theta_\rho(\hat\theta)$. This keeps the certificate statistically meaningful even for large, highly parameterized generators.

%When this uniform slack is strictly positive, part of the finite-simulation error is absorbed by the conservatism induced by robustification. When $\bar\mu(\rho;x)=0$, the theorem remains valid but reduces to the usual coverage-based certificate without slack improvement. Thus, the theorem should be interpreted in the following way: the presence of a true-distribution-equivalent sampler in $\Theta_\rho(\hat\theta)$ protects against sampler misspecification, while any positive uniform margin between the robust and true objectives tightens the simulation guarantee. The benefit from robustification therefore depends not only on coverage, but also on whether the induced slack is nontrivial uniformly over the decision class.

\subsection{Proof of Theorem \ref{thm:sro_generalization}}

Under (A0), for every $\omega\in\mathcal W$,
\[
J_{\tilde\theta^\star}(\omega;x)=J^\star(\omega;x),
\]
since $P_{\tilde\theta^\star}(\cdot\mid x)=P_{\theta^\star}(\cdot\mid x)$ by assumption. Because $\tilde\theta^\star\in\Theta_\rho(\hat\theta)$, it follows that
\[
U(\omega;x)
=
\sup_{\theta\in\Theta_\rho(\hat\theta)}J_\theta(\omega;x)
\ge
J_{\tilde\theta^\star}(\omega;x)
=
J^\star(\omega;x).
\]
Hence $\mu(\omega,\rho;x)\ge 0$ for all $\omega\in\mathcal W$, and therefore $\bar\mu(\rho;x)\ge 0$.

Fix $\omega\in\mathcal W$. By (A5), let
\[
\bar\theta(\omega)\in\Theta_\rho(\hat\theta)
\quad\text{satisfy}\quad
U(\omega;x)=J_{\bar\theta(\omega)}(\omega;x).
\]
Since $\bar\theta(\omega)$ is defined through population expectations, it is independent of the latent samples $\{z_i\}_{i=1}^N$. By Gaussian concentration for Lipschitz functions, for any $\delta'\in(0,1)$, with probability at least $1-\delta'$,
\begin{equation}
J_{\bar\theta(\omega)}(\omega;x)
\le
\frac{1}{N}\sum_{i=1}^N f(\omega,G_{\bar\theta(\omega)}(z_i,x))
+
L_yL_z\sqrt{\frac{2\log(1/\delta')}{N}}.
\label{eq:appendix_fixed_omega}
\end{equation}
Because $\bar\theta(\omega)\in\Theta_\rho(\hat\theta)$,
\begin{equation}
\frac{1}{N}\sum_{i=1}^N f(\omega,G_{\bar\theta(\omega)}(z_i,x))
\le
\hat U_N(\omega;x).
\label{eq:appendix_empirical_domination}
\end{equation}
Combining \eqref{eq:appendix_fixed_omega} and \eqref{eq:appendix_empirical_domination} yields
\begin{equation}
J_{\bar\theta(\omega)}(\omega;x)
\le
\hat U_N(\omega;x)
+
L_yL_z\sqrt{\frac{2\log(1/\delta')}{N}}.
\label{eq:appendix_worst_case}
\end{equation}
Since
\[
J_{\bar\theta(\omega)}(\omega;x)=U(\omega;x)=J^\star(\omega;x)+\mu(\omega,\rho;x),
\]
we obtain
\begin{equation}
J^\star(\omega;x)
\le
\hat U_N(\omega;x)
+
L_yL_z\sqrt{\frac{2\log(1/\delta')}{N}}
-
\mu(\omega,\rho;x).
\label{eq:appendix_pointwise}
\end{equation}

We now uniformize over $\omega\in\mathcal W$. Let $\mathcal M_\varepsilon\subseteq\mathcal W$ be an $\varepsilon$-cover of size $N_{\mathcal W}(\varepsilon)$. Apply \eqref{eq:appendix_pointwise} to each $\omega_j\in\mathcal M_\varepsilon$ with
\[
\delta'=\frac{\delta}{N_{\mathcal W}(\varepsilon)}.
\]
By a union bound, with probability at least $1-\delta$, simultaneously for all $\omega_j\in\mathcal M_\varepsilon$,
\begin{equation}
J^\star(\omega_j;x)
\le
\hat U_N(\omega_j;x)
+
L_yL_z
\sqrt{
\frac{
2\bigl(\log N_{\mathcal W}(\varepsilon)+\log(1/\delta)\bigr)
}{N}}
-
\mu(\omega_j,\rho;x).
\label{eq:appendix_net_bound}
\end{equation}

Now fix any $\omega\in\mathcal W$, and choose $\omega_j\in\mathcal M_\varepsilon$ such that $\|\omega-\omega_j\|\le\varepsilon$. By (A2),
\begin{equation}
|J^\star(\omega;x)-J^\star(\omega_j;x)|
\le
L_\omega\varepsilon.
\label{eq:appendix_population_transfer}
\end{equation}
Likewise, for every $\theta\in\Theta_\rho(\hat\theta)$,
\[
\left|
\frac{1}{N}\sum_{i=1}^N f(\omega,G_\theta(z_i,x))
-
\frac{1}{N}\sum_{i=1}^N f(\omega_j,G_\theta(z_i,x))
\right|
\le
L_\omega\varepsilon,
\]
and therefore
\begin{equation}
\hat U_N(\omega_j;x)
\le
\hat U_N(\omega;x)+L_\omega\varepsilon.
\label{eq:appendix_empirical_transfer}
\end{equation}
Finally, since $\mu(\omega_j,\rho;x)\ge \bar\mu(\rho;x)$,
\begin{equation}
-\mu(\omega_j,\rho;x)\le -\bar\mu(\rho;x).
\label{eq:appendix_slack_transfer}
\end{equation}
Substituting \eqref{eq:appendix_population_transfer}, \eqref{eq:appendix_empirical_transfer}, and \eqref{eq:appendix_slack_transfer} into \eqref{eq:appendix_net_bound} yields, simultaneously for all $\omega\in\mathcal W$,
\[
J^\star(\omega;x)
\le
\hat U_N(\omega;x)
+
L_yL_z
\sqrt{
\frac{
2\bigl(\log N_{\mathcal W}(\varepsilon)+\log(1/\delta)\bigr)
}{N}}
-
\bar\mu(\rho;x)
+
2L_\omega\varepsilon.
\]
This proves Theorem~\ref{thm:sro_generalization}.

\section{Additional Experimental Details}
\label{app:exp_details}
\subsection{Architectural Specifications}
\label{subsec:architecture}

\begin{figure}[H]
    \centering
    \includegraphics[width=0.8\textwidth]{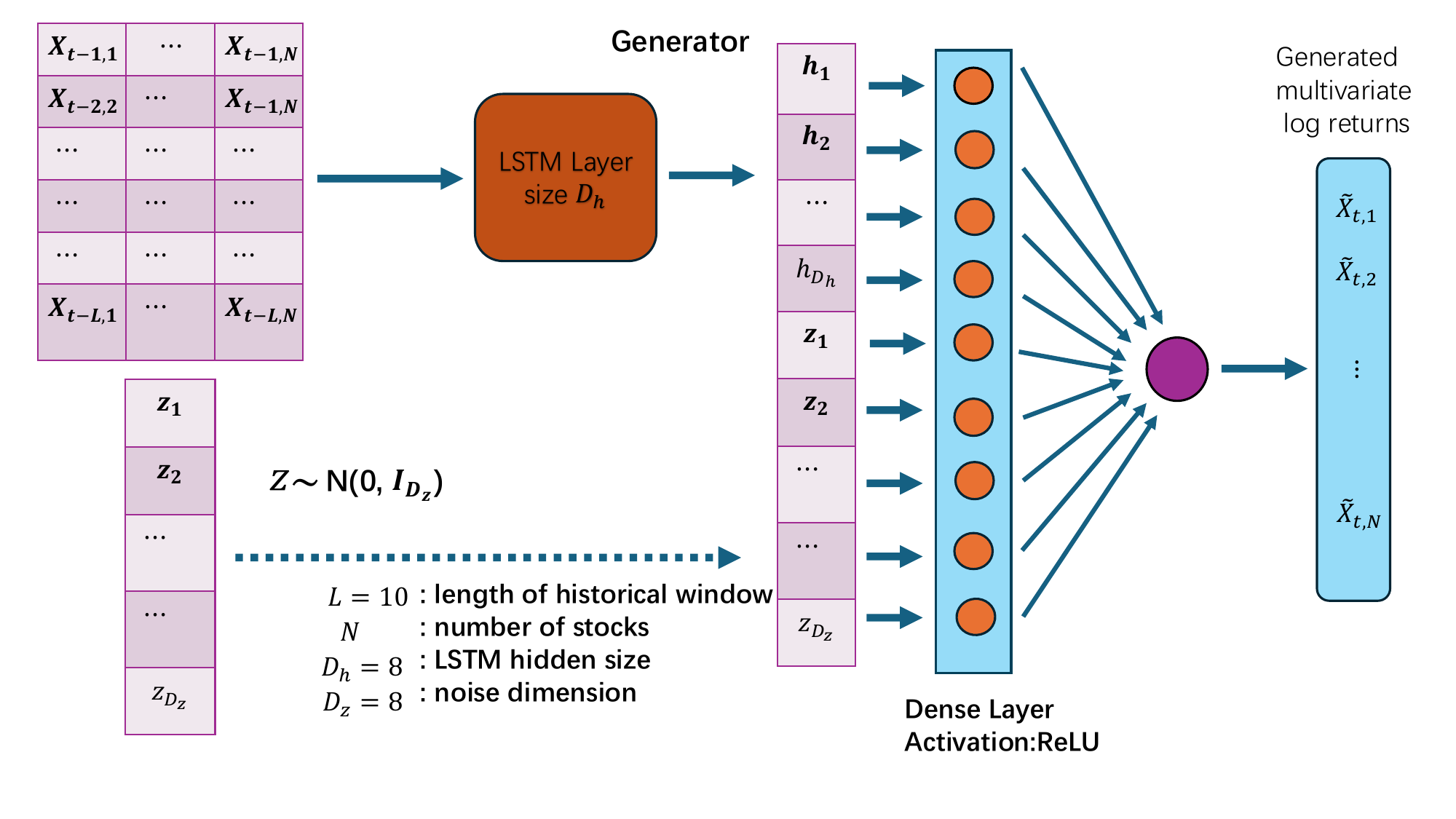}
    \caption{Architecture of the conditional LSTM-based generator.
    The generator takes a historical log-return matrix $X_{i,j} \in \mathbb{R}^{L \times N}$, 
where $i$ indexes time and $j$ indexes assets, together with a noise vector $z$, 
and outputs simulated next-day multivariate log returns. In the controlled experiment, $N=20$, whereas in the real-data experiment, $N=10$.}
    \label{fig:generator}
\end{figure}

\begin{figure}[H]
    \centering
    \includegraphics[width=0.8\textwidth]{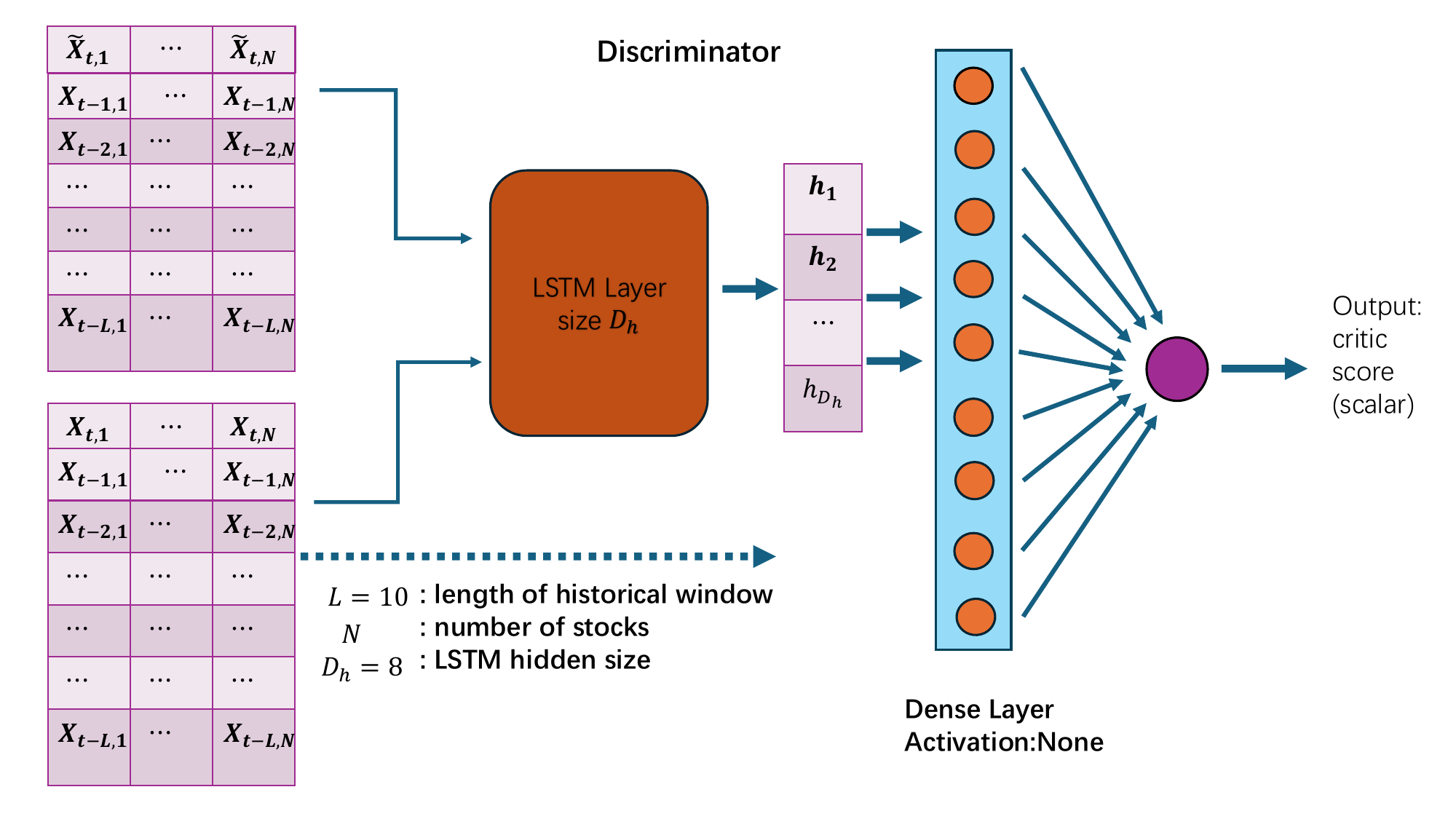}
    \caption{Architecture of the conditional LSTM-based discriminator (critic). 
The discriminator takes the historical log-return matrix $X_{i,j}$ as conditioning input and concatenates it 
with either the realized or generated next-day return vector. 
The resulting sequence is processed by an LSTM and mapped to an unconstrained scalar critic score 
for real-versus-synthetic discrimination.}
    \label{fig:discriminator}
\end{figure}

The nominal sampler $G_{\hat\theta}$ follows the recurrent topology of \cite{vuletic2024fin}, while replacing the original adversarial loss with a least-squares GAN objective for numerical stability and computational tractability. Accordingly, we remove the terminal sigmoid activation from the discriminator $D_\eta$ in Fig.~\ref{fig:discriminator} and use the corresponding squared-loss objective. The resulting architecture follows a standard encoder--decoder structure:

\paragraph{Contextual Encoding.} 
The encoder $E_\phi$ processes the conditioning window $\mathbf{x} \in \mathbb{R}^{L \times N}$ through an LSTM and maps it to a low-dimensional contextual representation:
\begin{equation}
        \hat{H} = E(\mathbf{x}; \phi) := \mathrm{LSTM}(\mathbf{x}_{1:L}; \phi),
\end{equation}
where $\mathbf{x}$ is standardized using the historical mean and standard deviation estimated from the first 100 observations of the training window, so that the normalization step does not use future information. The terminal hidden state $\hat{H}$ is used as a contextual representation of the conditioning window.

\paragraph{Stochastic Decoding.} 
As illustrated in Fig.~\ref{fig:generator}, the decoder $G_\psi$ combines the contextual representation $\hat{H}$ with latent noise $\mathbf{z} \sim \mathcal{N}(0, \mathbf{I})$ and maps the concatenated vector to a simulated next-period multivariate return:
\begin{equation}
    \tilde{\mathbf{y}} = \mathrm{ReLU}(W_\psi [\hat{H}; \mathbf{z}] + b_\psi) \odot s + m,
\end{equation}
where $m$ and $s$ denote the mean and standard deviation used in the normalization step. Thus, the decoder output is passed through a ReLU activation and then transformed back to the original return scale.
\subsection{Calibration Protocol and Hyperparameters}

The generator-side parameters $\hat{\theta}=\{\phi,\psi\}$ are calibrated on the training set $\mathcal{D}$ using a least-squares adversarial objective. Algorithm~\ref{alg:lsgan_calibration} summarizes the corresponding adversarial training procedure, in which the generator learns to produce realistic return scenarios while the discriminator provides the training signal through the least-squares objective. The main network dimensions and training hyperparameters, including the learning rates and latent dimensions, are summarized in Table~\ref{tab:nominal_gan_config}.

\begin{algorithm}[H]
\caption{Nominal Sampler Calibration (LSGAN)}
\label{alg:lsgan_calibration}
\small
\begin{algorithmic}[1]
\Require Training set $\mathcal{D} = \{(\mathbf{x}, \mathbf{y})\}$; noise primitive $\mathbf{z} \sim \mathcal{N}(0, \mathbf{I})$.
\Ensure Optimized sampler parameters $\hat{\theta} = \{\phi, \psi\}$.
\State Initialize sampler $G_\theta$ (with $\theta = \{\phi, \psi\}$) and auxiliary critic $D_\eta$.
\State Set calibration targets: $\tau_{r} \gets 1$, $\tau_{f} \gets 0$.
\For{epoch $1$ to $E$}
    \For{mini-batch $(\mathbf{x}, \mathbf{y}) \sim \mathcal{D}$}
        \State Sample noise $\mathbf{z} \sim \mathcal{N}(0, \mathbf{I})$ and generate $\tilde{\mathbf{y}} = G_\theta(\mathbf{z}, \mathbf{x})$.
        \State Update auxiliary critic $\eta$ by minimizing the quadratic loss:
        \Statex \hskip\algorithmicindent $\mathcal{L}_D(\eta) = \frac{1}{2} \mathbb{E} \left[ (D_\eta([\mathbf{x}; \mathbf{y}]) - \tau_{r})^2 + (D_\eta([\mathbf{x}; \tilde{\mathbf{y}}]) - \tau_{f})^2 \right]$
        \State Update sampler $\theta$ by minimizing the adversarial loss:
        \Statex \hskip\algorithmicindent $\mathcal{L}_G(\theta) = \frac{1}{2} \mathbb{E} \left[ (D_\eta([\mathbf{x}; G_\theta(\mathbf{z}, \mathbf{x})]) - \tau_{r})^2 \right]$
    \EndFor
\EndFor
\State \Return $\hat{\theta} \gets \{\phi, \psi\}$ \Comment{Only the calibrated sampler is retained for the SRO framework}
\end{algorithmic}
\end{algorithm}

\begin{table}[H]
\centering
\caption{Training hyperparameters for the nominal conditional LSTM-GAN}
\label{tab:nominal_gan_config}
\small
\setlength{\tabcolsep}{8pt}
\renewcommand{\arraystretch}{1.15}
\begin{tabular}{ll}
\toprule
\textbf{Hyperparameter} & \textbf{Value} \\
\midrule
Optimizer & RMSprop \\
Generator learning rate & $5\times 10^{-5}$ \\
Discriminator learning rate & $3\times 10^{-5}$ \\
Batch size & 40 \\
Epochs & 2000 \\
Look-back window ($L$) & 10 \\
Hidden dimension ($D_h$) & 8 \\
Noise dimension ($D_z$) & 8 \\
Weight initialization & Xavier normal \\
Asset dimension ($d$) & 20 in the controlled experiment; 10 in the real-data experiment \\
Normalization & Mean/std estimated from the first 100 training observations \\
\bottomrule
\end{tabular}
\end{table}

\begin{table}[H]
\centering
\caption{SRO and experiment-specific hyperparameters}
\label{tab:sro_exp_config}
\small
\setlength{\tabcolsep}{7pt}
\renewcommand{\arraystretch}{1.15}
\begin{tabular}{ll}
\toprule
\textbf{Hyperparameter} & \textbf{Value} \\
\midrule
\multicolumn{2}{c}{\textbf{Common SRO settings}} \\
\midrule
Robustness radius ($\rho$) & 0.3 \\
Utility parameter ($\lambda$) & 10 \\
Learning rate ($\theta$-update) & 0.001 \\
Learning rate ($\omega$-update) & 0.1 \\
Iteration budget & 12000 \\
%CVaR level ($\alpha$) & 0.05 \\
\midrule
\multicolumn{2}{c}{\textbf{Controlled generator-to-generator experiment}} \\
\midrule
Asset dimension ($d$) & 20 \\
Training window & 1000 \\
Synthetic path length & 1000 \\
Retraining window & 800 \\
Validation window & 100 \\
Out-of-sample test window & 100 \\
Seeds & 40--59 \\
\midrule
\multicolumn{2}{c}{\textbf{Real-data experiment}} \\
\midrule
Asset dimension ($d$) & 10 \\
Training window & 1000 \\
Validation window & 100 \\
Out-of-sample test window & 100 \\
Seeds & 40--59 \\
Excluded run(s) & seed 49 \\
\bottomrule
\end{tabular}
\end{table}

\subsection{Stock universe and data preprocessing}
\label{subsec:stock_universe}

Our empirical study is based on a panel of 50 large-cap U.S. equities:
AAPL, MSFT, GOOGL, AMZN, META, TSLA, NVDA, ADBE, INTC, CRM, AMD, CSCO, ORCL, IBM, QCOM, JPM, BAC, WFC, C, GS, MS, V, MA, AXP, BLK, JNJ, PFE, MRK, UNH, ABBV, TMO, ABT, LLY, DHR, BMY, KO, PG, PEP, WMT, DIS, HD, MCD, NKE, SBUX, COST, XOM, CVX, VZ, T, and NFLX. The sample period runs from May 15, 2020 to May 14, 2025.

%For the controlled experiment, we fix a 20-stock universe consisting of AAPL, MSFT, GOOGL, AMZN, META, TSLA, NVDA, ADBE, INTC, CRM, AMD, CSCO, ORCL, IBM, QCOM, JPM, BAC, WFC, C, and GS. For the real-data experiment, we randomly sample 10-stock sub-universes from the 50-stock panel in each seeded run in order to reduce universe-specific dependence and maintain stable GAN training.

Daily prices are converted into log returns before model training. To reduce the influence of extreme observations and stabilize generator calibration, each log return is clipped to the interval $[-0.1,\,0.1]$. The resulting return panel is then used to construct the rolling conditioning windows and next-period targets for both the nominal generator and the downstream portfolio experiments.

\bibliography{references}
\end{document}